\documentclass[12pt]{amsart}
\usepackage{amsfonts, amssymb, amsmath, amsthm, latexsym, array}
\usepackage{fullpage,color}
\usepackage{verbatim,euscript}

\input {cyracc.def}
 \font\tencyr=wncyr10 
\font\tencyi=wncyi10 
\font\tencysc=wncysc10 
\def\rus{\tencyr\cyracc}
\def\rusi{\tencyi\cyracc}
\def\rusc{\tencysc\cyracc}

\newtheorem{thm}{Theorem}[section]

\newtheorem{lm}[thm]{Lemma}
\newtheorem{prop}[thm]{Proposition}

\theoremstyle{remark}
\newtheorem{ex}[thm]{Example}

\theoremstyle{definition}
\newtheorem{rmk}[thm]{Remark}
\newtheorem{df}{Definition}

\newcommand {\g}{{\mathfrak g}}

\newcommand {\q}{{\mathfrak q}}

\newcommand {\z}{{\mathfrak z}}

\newcommand {\glv}{{\mathfrak{gl}(V)}}
\newcommand {\gln}{{\mathfrak{gl}_n}}
\newcommand {\sln}{{\mathfrak{sl}_n}}

\newcommand {\slno}{{\mathfrak{sl}_{n+1}}}

\newcommand {\spn}{{\mathfrak {sp}}_{2n}}
\newcommand {\son}{\mathfrak{so}_{n}}

\newcommand {\ca}{{\mathcal A}}

\newcommand {\cf}{{\mathcal F}}
\newcommand {\ck}{{\mathcal K}}

\newcommand {\cs}{{\mathcal S}}

\newcommand {\BP}{{\mathbb P}}

\newcommand {\BZ}{{\mathbb Z}}

\newcommand {\md}{/\!\!/}

\newcommand{\lb}{\lambda}
\newcommand{\ap}{\alpha}

\renewcommand{\le}{\leqslant}
\renewcommand{\ge}{\geqslant}

\newcommand{\eus}{\EuScript}
\newcommand {\zq}{{\eus Z(\q)}}
\newcommand {\zg}{{\eus Z(\g)}}

\newcommand {\ad}{{\mathrm{ad\,}}}
\newcommand {\ads}{{\mathrm{ad}^*}}

\newcommand {\codim}{{\mathrm{codim\,}}}

\newcommand {\ind}{{\mathrm{ind\,}}}

\newcommand {\Mor}{\operatorname{Mor}}

\newcommand {\rk}{{\mathrm{rk\,}}}

\newcommand {\spe}{{\mathrm{Spec\,}}}

\newcommand {\trdeg}{{\mathrm{trdeg\,}}}

\newcommand {\GR}[2]{{\textrm{{\bf #1}}}_{#2}}

\newfam\eusfam%
\font\Bbbfont=msbm10 scaled 1200%
\def\bbk{\hbox {\Bbbfont\char'174}}

\begin{document}
\setlength{\parskip}{3pt plus 5pt minus 0pt}
\hfill {\scriptsize February 19, 2007}
\vskip1ex

\title[On maximal commutative subalgebras]
{The argument shift method and 
maximal commutative subalgebras of Poisson algebras}
\author[D.\,Panyushev]{Dmitri I.~Panyushev}
\address[D.P.]{Independent University of Moscow,
Bol'shoi Vlasevskii per. 11, 119002 Moscow, \ Russia}
\email{panyush@mccme.ru}
{\color{cyan}\author[O.\,Yakimova]{Oksana S.~Yakimova}}
\address[O.Y.]{Universit\"at zu K\"oln,
Mathematisches Institut, Weyertal 86-90, D-50931 K\"oln, Deutschland}
\email{oyakimov@math.uni-koeln.de}
\thanks{Both authors are supported by  R.F.B.R. grant 05--01--00988.
The second author  gratefully acknowledges the support of the
Alexander von Humboldt-Stiftung.}
\maketitle

\section*{Introduction}

\noindent
Let $\q$ be a Lie algebra over an algebraically closed field $\bbk$ of characteristic zero.
The symmetric algebra $\cs(\q)$ has a natural structure of Poisson algebra, and 
our goal is to present a sufficient condition for the maximality 
of  Poisson-commutative subalgebras of $\cs(\q)$ obtained by the
argument shift method. Study of Poisson-commuttive subalgebras of $\cs(\q)$
has attracted much attention in the last years, see \cite{jp, kw1, rybn, vitya, tar1}. 
This is related to commutative subalgebras of the enveloping algebra $\eus U(\q)$, 
fine questions of  
symplectic geometry, and
integrable Hamiltonian systems.
Commutative subalgebras of  $\eus U(\q)$ (e.g., the famous Gelfand-Zetlin
subalgebra of $\eus U(\sln)$)
occur in the theory of quantum integrable systems
and have interesting application in representation theory. 
 
Let  $\zq$ be the centre of the Poisson algebra $\cs(\q)$.
For $\xi\in\q^*$, let $\cf_\xi(\zq)$ denote the
algebra generated by the $\xi$-shifts of all $f\in\zq$ (see Subsection~\ref{sdvig} for precise
definitions). As is well-known, $\cf_\xi(\zq)$ is a Poisson-commutative subalgebra of
$\cs(\q)$. Furthermore, $\trdeg(\cf_\xi(\zq))\le(\dim\q+\ind\q)/2=:b(\q)$.
We say that $\cf_\xi(\zq)$ is  {\it of maximal dimension}, if the equality holds.
However, even in this case, it may happen that there is a strictly larger 
Poisson-commutative subalgebra (of the same transcendence degree).
We say that $\cf_\xi(\zq)$ is {\it maximal}, if it is 
maximal with respect to inclusion among the commutative subalgebras of $\cs(\q)$.
Let $\q^*_{reg}$ denote the set of {\it regular\/} elements of $\q^*$, i.e., those whose 
stabiliser in $\q$ has the minimal dimension.
For the purposes of this introduction, we state our main result 
(Theorem~\ref{thm:cod3}) in a slightly abbreviated form:

\begin{thm}  \label{intro_main}
Suppose that
\begin{itemize}
\item[\sf (i)] \ 
$\eus Z(\q)$ contains algebraically independent homogeneous polynomials
$f_1,\ldots,f_l$, where $l=\ind\q$, such that $\sum_{i=1}^l \deg f_i=b(\q)$;
\item[\sf (ii)] \ $\codim (\q^*\setminus \q^*_{reg})\ge 3$.
\end{itemize}
Then, for any $\xi\in \q^*_{reg}$,  $\cf_\xi(\eus Z(\q))$ is a  polynomial algebra
of Krull dimension $b(\q)$ and it is a maximal Poisson-commutative subalgebra of 
$\cs(\q)$.
\end{thm}
\noindent
Obviously, Theorem~\ref{intro_main} applies if $\q$ is semisimple, and we thus 
generalise results 
of A.\,Tarasov \cite{tar1}. (He proved maximality if $\xi$ is regular semisimple.)
There are also other interesting classes of Lie algebras satisfying the conditions
of this theorem, see Section~\ref{sect:applic}.

A general criterion of Bolsinov~\cite{bol2} asserts that, 
for  $\xi\in \q^*_{reg}$, \ $\cf_\xi(\eus Z(\q))$
is of maximal dimension if and only if $\codim (\q^*\setminus \q^*_{reg})\ge 2$.
For the proof of Theorem~\ref{intro_main}, we need, however, a stronger 
result. Namely, we provide a precise
description of pairs $\xi,\eta\in\q^*$
such that the differentials at $\eta$ of all functions from $\cf_\xi(\eus Z(\q))$
generate a subspace of dimension $b(\q)$, see Theorem~\ref{compl}.

{\sl Notation.}
If an algebraic group $Q$ acts on an irreducible affine variety $X$, then $\bbk[X]^Q$ 
is the algebra of $Q$-invariant regular functions on $X$ and $\bbk(X)^Q$
is the field of $Q$-invariant rational functions. If $\bbk[X]^Q$
is finitely generated, then $X\md Q:=\spe \bbk[X]^Q$, and
the {\it quotient morphism\/} $\pi_X: X\to X\md Q$ is the mapping associated with
the embedding $\bbk[X]^Q \hookrightarrow \bbk[X]$.
\\   \indent
If $V$ is a $Q$-module and $v\in V$, then $\q_v$ is the stabiliser of 
$v$ in $\q$.
For the adjoint representation of $\q$, the stabiliser of $x\in \q$ is also denoted by 
$\z_\q(x)$, and we say that $\z_\q(x)$ is the {\it centraliser\/} of $x$.
\\   \indent
All topological terms refer to the Zariski topology.
If $M$ is a subset of a vector space, then ${\rm span}(M)$ denotes the linear span of
$M$; \ $\bbk^\times:=\bbk\setminus\{0\}$.

\section{On the codim--$n$ property for the coadjoint representation} 
\label{prelim}

Let $Q$ be a connected algebraic group with Lie algebra $\q$.
We write $\cs(\q)$ for the symmetric algebra of $\q$.
Recall that $\mathcal S(\q)\simeq \bbk[\q^*]$ is a Poisson algebra, and the symplectic leaves
in $\q^*$ are precisely the coadjoint orbits of $Q$. Since each coadjoint
orbit $Q{\cdot}\xi$ is a symplectic variety, $\dim Q{\cdot}\xi$ is even.
 Let $\{\ ,\ \}$ denote the Lie-Poisson
bracket in $\cs(\q)$. Then the algebra of invariants
$\bbk[\q^*]^Q=\cs(\q)^Q$ is the centre of $(\cs(\q), \{\ ,\ \})$.
We also write $\eus Z(\q)$ for this centre.

Let $\q^*_{reg}$ denote the set of all $Q$-{\it regular\/} elements of $\q^*$. That is,
\[
   \q^*_{reg}=\{\xi\in\q^*\mid \dim Q{\cdot}\xi\ge \dim Q{\cdot}\eta 
   \text{ for all } \eta\in\q^*\} \ .
\]
As is well-known, $\q^*_{reg}$ is a dense open subset of $\q^*$.

\begin{df}   \label{def:codim2}
We say that the coadjoint representation of $\q$ has the {\it codim--n property\/} if
$\codim (\q^*\setminus \q^*_{reg})\ge n$.
\end{df}

\noindent
If $\xi\in \q^*_{reg}$, then $\dim\q_\xi$ is called the {\it index\/} of $\q$, denoted $\ind\q$.
By Rosenlicht's theorem, $\trdeg\bbk(\q^*)^Q=\ind\q$. It follows that
if $f_1,\dots,f_r\in \bbk[\q^*]^Q$ are algebraically independent, then $r\le \ind\q$.
Set $b(\q)=(\dim\q+\ind\q)/2$.
If $\q$ is semisimple, then $b(\q)$ is the dimension of  a Borel subagebra.

\noindent
{\it  Example.} If $\g$ is reductive, then $\ad\simeq\ads$ and 
$\codim (\g\setminus \g_{reg})=3$. Hence the coadjoint representation of a 
reductive Lie algebra has the codim--$3$ property.

The following example pointed out by E.B.\,Vinberg
shows that for any $n$ there are noncommutative Lie algebras with codim--$n$ property.

\begin{ex}    \label{ex:vinberg}
Suppose $s\in \glv$ is a semisimple linear transformation
with nonzero rational eigenvalues.
Let $\q$ be the semi-direct product of the 1-dimensional toral Lie algebra $\bbk s$ and $V$.
The Lie bracket is given by
\[
     [(\ap s,v),(\beta s,v')]= (0, \ap s(v')- \beta s(v)), \qquad \ap,\beta\in\bbk.
\]
It is easily seen that $\ind \q=\dim\q -2$. Moreover, let $L$ be the annihilator of $V$ in $\q^*$.
Then the line $L$ is precisely the set of $Q$-fixed points in $\q^*$, while 
$\dim Q{\cdot}\xi=2$ for any $\xi\in \q^*\setminus  L$. Thus, $\q$ has the codim--$n$
property with $n=\dim V$.
\end{ex}

If $f\in \cs(\q)$, then the differential of $f$, $\textsl{d}f$, can be regarded as a  
polynomial mapping from $\q^*$ to $\q$, i.e., an element of
$\Mor_Q(\q^*,\q)\simeq \cs(\q)\otimes \q$. 
More precisely, if 
$f\in \cs^d(\q)$, then
$\textsl{d}f$ is a polynomial mapping of degree $d-1$, i.e., an element
of $\cs^{d-1}(\q)\otimes \q$.
We write $(\textsl{d}f)_\xi$ for the value of $\textsl{d}f$ at $\xi\in\q^*$.
Recall that $(\textsl{d}f)_\xi$ is an element of $\q$
that is defined as follows. If $\nu\in \q^*$ and $\langle\ ,\ \rangle$ denotes the natural pairing
between $\q$ and $\q^*$, then 
\[
\langle (\textsl{d}f)_\xi,\nu\rangle:= \text{the coefficient of
$t$ in the Taylor expansion of $f(\xi+t\nu)$}.
\] 
The r\^ole of the codim--$2$ property is seen in the following result,
see \cite[Theorem\,1.2]{coadj}.

\begin{thm}  \label{thm:cod2}
Suppose that $(\q,\ads)$ has the codim--2 property
and $\trdeg \bbk[\q^*]^Q=\ind\q$. Set  $l=\ind\q$. 
Let $f_1,\dots,f_l\in \bbk[\q^*]^Q$ be arbitrary
homogeneous algebraically independent polynomials. Then 
\begin{itemize}
\item[\sf (i)] \ $\sum_{i=1}^l\deg f_i \ge b(\q)$; 
\item[\sf (ii)] \  If\/  $\sum_{i=1}^l\deg f_i = b(\q)$, then
$\bbk[\q^*]^Q$ is freely generated by $f_1,\dots,f_l$ and
$\xi\in\q^*_{reg}$ if and only if $(\textsl{d}f_1)_\xi,\dots,(\textsl{d}f_l)_\xi$ are linearly
independent.
\end{itemize}
\end{thm}%
\noindent
The second assertion in (ii) can be regarded as a generalisation of Kostant's result for
reductive Lie algebras \cite[(4.8.2)]{ko63}.
Its geometric meaning is the following. Consider
the quotient morphism $\pi: \q^*\to \q^*\md Q\simeq \mathbb A^{\ind\q}$.
Then $\pi$ is smooth at $\xi\in\q^*$ if and only if $\xi\in\q^*_{reg}$.

\section{The argument shift  method and Bolsinov's criterion}
\label{sect:shift}

\subsection{Commutative subalgebras of $\cs(\q)$}
Let $\eus A$ be a subalgebra of the symmetric algebra $\cs(\q)$.
Then $\eus A$ is said to be {\it Poisson--commutative} if the restriction of $\{\ ,\ \}$ to 
$\eus A$ is zero. Abusing the language, we will usually omit "Poisson"
and merely say that $\eus A$ is commutative.
Notice that the words "subalgebra of $\cs(\q)$" always refer to the usual 
(associative and commutative) structure of the symmetric
algebra, while  "commutative" refers to the Poisson structure on $\cs(\q)$.

For any subalgebra $\eus A\subset \cs(\q)$, we define the transcendence degree
of $\eus A$ as that of the quotient field of $\eus A$.
As is well-known, if $\eus A$ is commutative, then
$\trdeg \eus A\le b(\q)$. Indeed, if $f_1,\dots, f_n\in \eus A$ are algebraically independent, 
then $M:={\rm span}\{(\textsl{d}f_1)_\xi,\dots,(\textsl{d}f_n)_\xi\}$ 
is $n$-dimensional for generic $\xi$. Furthermore, $M$ is an isotropic subspace of
$\q$ with respect to the Kirillov form $\ck_\xi$.
(Recall that $\ck_\xi(x,y):=\langle \xi, [x,y]\rangle$ and hence $\dim(\ker\ck_\xi)=\dim\q_\xi$.)

\begin{df} Let $\ca$ be a commutative subalgebra of $\cs(\q)$. 
Then $\ca$ is said to be {\it of maximal dimension}, if
$\trdeg\ca=b(\q)$; \ $\ca$ is said to be {\it maximal}, if 
it is maximal with respect to inclusion among the commutative subalgebras of $\cs(\q)$.
\end{df}

\noindent We do not know of whether there exist maximal commutative subalgebras
that are not of maximal dimension.
\\  \indent
Suppose $\ca$ is commutative and of maximal dimension.
If $\ca\subset\ca'$ and $\ca'$ is commutative, then each element of $\ca'$ is algebraic
over $\ca$. Conversely, if  $f\in\cs(\q)$ is algebraic over $\ca$, then, for generic
$\xi\in\q^*$, $(\textsl{d}f)_\xi$ belongs to ${\rm span}\{(\textsl{d}F)_\xi \mid F\in \ca\}$,
which is an isotropic subspace with resepect to $\ck_\xi$. Hence 
$\{f,F\}(\xi)=0$ for a generic $\xi$ and therefore $\{f,F\}\equiv 0$.
Thus, $\ca$ is maximal if and only if it is algebraically closed in $\cs(\q)$.

\subsection{The argument shift method}      \label{sdvig}
Suppose $f\in\cs(\q)$ is a  polynomial of degree $d$.
For any $\xi\in\q^*$, we may consider a shift of $f$ in direction $\xi$:
$f_{a,\xi}(\mu)=f(\mu+a\xi)$, where $a\in \bbk$.
Expanding the right hand side as polynomial in $a$, we obtain the expression
$f_{a,\xi}(\mu)=\sum_{j=0}^d f_{\xi}^j(\mu)a^j$.
Associated with this shift of argument, we obtain the family of polynomials
$ f_{\xi}^j$, where $j=0,1,\dots,d-1$. (Since $\deg f_{\xi}^j=d-j$, the value
$j=d$ is not needed.)
We will say that the polynomials $\{ f_{\xi}^j\}$ are $\xi$-{\it shifts\/} of $f$.
Notice that $f_\xi^0=f$ and $f_\xi^{d-1}$ is a linear form on $\q^*$, i.e.,
an element of $\q$. Actually, $f_\xi^{d-1}=(\textsl{d}f)_\xi$.
There is also an obvious symmetry with respect to $\xi$ and $\mu$: \ 
$f_\xi^j(\mu)=f_\mu^{d-j}(\xi)$.

The following observation is due to Mishchenko-Fomenko \cite{mf3}.

\begin{lm}     \label{mf:shift}
Suppose that $h_1,\dots, h_m\in \eus Z(\q)$. 
Then for any $\xi\in\q^*$, the polynomials
\[
  \{h_{i,\xi}^j \mid i=1,\dots,m;\quad j=0,1,\dots,\deg h_i-1\} 
\]
pairwise commute with respect to the Poisson bracket.
\end{lm}

\noindent
Mishchenko and Fomenko used this procedure for constructing
commutative subalgebras of maximal dimension in $\cs(\q)$.
Given $\xi\in\q^*$ and an arbitrary subset $\eus B\subset \zq$,
let  $\cf_\xi(\eus B)$ denote the subalgebra of $\cs(\q)$ generated 
by the $\xi$-shifts of all elements of $\eus B$. Clearly, if $\hat{\eus B}$ is the subalgebra
generated by $\eus B$, then $\cf_\xi(\eus B)=\cf_\xi(\hat{\eus B})$.
By Lemma~\ref{mf:shift}, all  subalgebras $\cf_\xi(\eus B)$ are  commutative.
In particular, subalgebras $\cf_\xi(\zq)$   are natural candidates on
the r\^ole of commutative subalgebras of maximal dimension.

For $\g$ semisimple, it is proved in \cite{mf3} that  there is an open subset
$\Omega\subset\g^*$ such that  $\cf_\xi(\zg)$ is of 
maximal dimension for any $\xi\in\Omega$.
The subalgebras of the form  $\cf_\xi(\zg)$ are called {\it Mishchenko-Fomenko subalgebras}
in \cite{tar1, vi90, vitya}.

\begin{rmk}   \label{2-skobki}
The argument shift method is a particular case of a more general construction related 
to compatible Poisson brackets. Recall that two Poisson brackets 
on a commutative associative algebra $\cs$ are said 
to be {\it compatible\/} if any linear combination of them is again a Poisson bracket.
For $\cs=\cs(\q)$, we can
consider the usual Lie-Poisson bracket
$(f,g)\to \{f,g\}$ and the bracket $(f,g)\to \{f,g\}_\xi$
obtained by  ``freezing the argument''. Here $f,g\in\cs(\q)$ and
$\xi\in\q^*$ is a fixed element. By definition, 
$\{f,g\}(\eta):=\langle \eta, [(\textsl{d}f)_\eta, (\textsl{d}g)_\eta]\rangle$ 
and $\{f,g\}_\xi(\eta):=\langle \xi, [(\textsl{d}f)_\eta, (\textsl{d}g)_\eta]\rangle$.
A direct calculation shows that 
each linear combination $a\{\,,\,\}+b\{\,,\,\}_\xi$ is again a Poisson bracket 
on ${\mathcal S}(\q)$. 

It is easily seen that if $f\in\zq$ and $f_{b,\xi}(\nu):=f(\nu+b\xi)$, then $f_{b,\xi}$
is a central function with respect to $\{\,,\,\}+b\{\,,\,\}_\xi$. Furthermore, the
assignment $f\mapsto f_{b,\xi}$ is a bijection between two centres. It follows that 
$\cf_\xi(\zq)$ is the subalgebra of $\cs(\q)$ generated by the centres of all
Poisson brackets $\{\,,\,\}+b\{\,,\,\}_\xi$, $b\in\bbk$.
\end{rmk}

\subsection{Bolsinov's criterion and its extension}  
A general criterion for $\cf_\xi(\zq)$ to be
of maximal dimension is found by A.V.~Bolsinov. 
Using our terminology, we can express it as follows. 

\begin{thm}[cf. Bolsinov {\cite[Theorem~3.1]{bol2}}]   \label{thm:bols}
Suppose that $\q$ satisfies the codim--2 property and $\trdeg\zq=\ind\q$. Then 
the algebra 
$\cf_\xi(\zq)$ is of maximal dimension for any $\xi\in \q^*_{reg}$.
\end{thm}

\noindent  The algebra $\cf_\xi(\zq)$ is of maximal dimension  if and only if
there is an  
$\eta\in\q^*$ such that the differentials at $\eta$ of all polynomials in $\cf_\xi(\zq)$
span a subspace of dimension $b(\q)$. Clearly, such $\eta$ form an open subset of $\q^*$.
For our main result,   we need, however, a more precise assertion.
Here it is.

\begin{thm}    \label{compl}
Keep the assumptions of Theorem~\ref{thm:bols}.
Let  $P\subset \q^*$ be a plane such that  
$P\setminus\{0\}\subset \q^*_{reg}$. Suppose that   
\\
\hbox to \textwidth{
\ $(\ast)$ \hfil  
$\dim{\rm span}\{(\textsl{d}f)_{\xi_0}\mid f\in\zq\}=\ind\q$ \ for some $\xi_0\in P$. \hfil }
Then
$\dim{\rm span}\{(\textsl{d}f)_{\eta} \mid f\in \cf_\xi(\zq)\}=b(\q)$
 for any linearly independent $\xi$ and $\eta$ in $P$.
\end{thm}

\noindent
{\it Remark.} Condition $(\ast)$ is open, hence it is satisfied on 
an open subset of $P$. In many important cases,  this condition 
follows from the other ones (see below). Therefore, there is not much harm in it.

\begin{proof}
We apply results of Bolsinov~\cite{bol2} (presented in Appendix~\ref{app}) 
to the compatible Poisson brackets 
$\{\,,\,\}$ and $\{\,,\,\}_\xi$ on $\q^*$, cf. Remark~\ref{2-skobki}.
For $\eta\in\q^*$, let 
$A_\eta$ and $B_\eta$ be the corresponding skew-symmetric forms on 
$T^*_\eta(\q^*)\cong\q$. 
Explicitly, $A_\eta(x,y)=\langle \eta,[x,y]\rangle$ and 
$B_\eta(x,y)=\langle \xi,[x,y]\rangle $. 
It follows that $(aA_\eta+b B_\eta)(x,y)=\langle a\eta+b\xi, [x,y]\rangle$ and
hence 
\begin{equation}  \label{ravno}
\dim (\ker (aA_\eta+b B_\eta))=\dim\q_{a\eta+b\xi} \ .
\end{equation}
We will  identify the $2$-dimensional vector spaces
$\eus P={\rm span}\{A_\eta,B_\eta\}$ and  $P={\rm span}\{\eta,\xi\}\subset \q^*$
by taking $a A_\eta+b B_\eta$ to $a\eta+b\xi$.

Set $\eus D:={\rm span}\{(\textsl{d}f)_\eta \mid f\in\cf_\xi(\zq)\}$. Our goal is
to prove that $\dim \eus D=b(\q)$.
Recall that $\trdeg \cs(\q)^{Q}=\ind\q$. Therefore
\[
   \Omega:=\{\nu\in \q^* \mid
   \dim{\rm span}\{(\textsl{d}f)_\nu\mid f\in {\mathcal S}(\q)^{Q} \} =\ind\q\}
\]
is a non-empty open subset of $\q^*$. Note that $\Omega$ is {\it conical}, i.e.,
$\nu\in \Omega$ if and only if $t\nu\in\Omega$ for any $t\in\bbk^\times$.
By the assumption, $\Omega_P:=\Omega\cap  P\ne \varnothing$.
\\
From Eq.~\eqref{ravno}, it follows that all nonzero forms in $\eus P$ have the same rank.
Applying Proposition~\ref{com1} to $V=\q$ and  
$\eus P={\rm span}\{A_\eta,B_\eta\}$ shows that  
$L=\sum_{(a,b)\ne(0,0)} \ker (aA_\eta+bB_\eta)$
is a maximal isotropic subspace of $\q$ with respect to any nonzero element of $\eus P$.
In particular, $\dim L=b(\q)$.
Furthermore, since $\Omega_P$ is a non-empty and conical subset
of $P\setminus\{0\}$, we deduce from
Lemma~\ref{open} that
\begin{equation}       \label{L}
L=\sum_{(1,b)\in \Omega_P} \ker (A_\eta+b B_\eta),
\end{equation}
 where 
$(1,b)$ is regarded as the point $\eta+b\xi\in P$.
Because $\dim \eus D\le b(\q)$, it 
suffices to prove  that $L\subset \eus D$. 
Take any $(1,b)\in \Omega_P$ and let $C=\{\,,\,\}+b\{\,,\,\}_\xi$ be 
the corresponding Poisson bracket on $\q^*$.
For any $f\in\zq$, set $\tilde f(\nu):=f(\nu+b\xi)$. Then 
$(\textsl{d}\tilde f)_\eta = (\textsl{d}f)_{\eta+b\xi}$ and 
$f\mapsto \tilde f$ is a 
bijection between $\zq$ and $\eus Z_C(\q)$, the centre of the Poisson algebra $(\cs(\q), C)$.
Hence 
\[   
\eus H:={\rm span}\{(\textsl{d}f)_{\eta+b\xi} \mid f\in\zq\}=
{\rm span}\{(\textsl{d}\tilde f)_\eta \mid f\in\eus Z_C(\q)\}\subset \ker(A_\eta+b B_\eta).
\]
Since $\eta+b\xi\in \Omega_P$, we have $\dim \eus H=\ind\q=\dim(\ker(A_\eta+b B_\eta))$. 
Hence 
${\rm span}\{(\textsl{d}\tilde f)_\eta \mid f\in\eus Z_C(\q)\}= \ker(A_\eta+b B_\eta)$.
But each $(\textsl{d}\tilde f)_\eta$ is a linear combination of differentials of elements of 
$\cf_\xi$. Therefore $\ker(A_\eta+b B_\eta)\subset \eus D$ whenever $(1,b)\in\Omega_P$, and we conclude from  Eq.~\eqref{L} that $L\subset \eus D$.  
Hence $L=\eus D$, and we are done.
\end{proof}

\section{Maximal commutative subalgebras of $\cs(\q)$ and flatness} 
\label{sect:max}

\noindent
First, we prove an auxiliary geometric result.
Let $V$ be a finite-dimensional vector space and $P\subset V$ a plane.
Suppose $\Omega$ is a conical open subset
of $V\setminus\{0\}$ such that $\codim(V\setminus\Omega)\ge n\ge 2$.
 Let us say that $P$ is an $\Omega$-plane
if $P\setminus \{0\}\subset \Omega$.
Given $v\in \Omega$, let $\Omega_v$ be the set of all $u$ such
that  $\bbk v+\bbk u\subset V$ is an $\Omega$-plane.

\begin{lm}   \label{auxil}
$\Omega_v$ is an open subset of\/ $V\setminus\{0\}$ 
and $\codim(V\setminus \Omega_v)\ge n-1$.
\end{lm}\begin{proof}
Set $S=V\setminus \Omega$ and consider 
the projectivisations $\BP(S)\subset \BP(V)$. Here $\BP(S)$ is a projective variety
of codimension $\ge n$. Write $\bar v$ for the image of  $v$ in $\BP(V)$.
Let $C$ be the cone in $\BP(V)$ generated by $\bar v$ and  $\BP(S)$.
That is, $C$ is the union of all lines through $\bar v$ and $y$, where $y$ runs over 
$\BP(S)$.
Then $C$ is a projective variety of codimension $\ge n-1$, and it follows from the 
construction that if $\bar u\not\in C$, then $\bbk v+\bbk u$ is an $\Omega$-plane.
Thus, $\BP(\Omega_v)=\BP(V)\setminus C$.
\end{proof}

The following is our main result.

\begin{thm}               \label{thm:cod3}
Let $\q$ be an algebraic Lie algebra.
\begin{itemize}
\item[\sf (i)] \ Suppose $(\q,\ads)$ has the codim--$2$ property
and $\eus Z(\q)$ contains algebraically independent polynomials
$f_1,\ldots,f_l$, where $l=\ind\q$, such that $\sum_{i=1}^l \deg f_i=b(\q)$.
Then, for any $\xi\in \q^*_{reg}$,  
$\cf_\xi(\eus Z(\q))=\cf_\xi(f_1,\ldots,f_l)$ is a  polynomial algebra
of Krull dimension $b(\q)$;

\item[\sf (ii)] \ Furthermore, if $(\q,\ads)$ has the codim--3 property, then 
$\cf_\xi(\eus Z(\q))$ is a maximal commutative subalgebra of $\cs(\q)$.
\end{itemize}
\end{thm}\begin{proof}
To simplify notation, write $\cf_\xi$ in place of $\cf_\xi(\eus Z(\q))$.

(i) It follows from the assumptions and Theorem~\ref{thm:cod2} that
$\zq=\bbk[f_1,\dots,f_l]$. Hence $\cf_\xi=\cf_\xi(f_1,\ldots,f_l)$.
By  Bolsinov's criterion (Theorem ~\ref{thm:bols}),  $\trdeg \mathcal F_\xi=b(\q)$ for any 
$\xi\in \q^*_{reg}$. 
Set $\Omega=\{\xi\in\q^*\mid (\textsl{d}f_1)_\xi,\dots, (\textsl{d}f_l)_\xi \text{ \ are linearly
independent}\}$. From Theorem~\ref{thm:cod2}(ii), it follows  that
$\Omega =\q^*_{reg}$. Hence $\codim (\q^*\setminus \Omega)\ge 2$.

Let $P:=\bbk\xi+\bbk\eta\subset \q^*$ be a $\q^*_{reg}$-plane, i.e., each nonzero element 
of it belongs to $\q^*_{reg}$. Since $\Omega=\q^*_{reg}$, each nonzero point of
$P$ satisfies condition $(\ast)$ of 
Theorem~\ref{compl}. Hence Theorem~\ref{compl} guarantees us that,
for any $\eta\in P\setminus \bbk\xi$,
the differentials of the $\xi$-shifts of $f_1,\ldots,f_l$
at $\eta$ span a subspace of dimension $b(\q)$.
Next, in view of the equality $\sum_{i=1}^l \deg f_i=b(\q)$, 
the set of all $\xi$-shifts of the $f_i$'s consists of
$b(\q)$ elements. It follows that  the differentials 
\[
   \{ (\textsl{d}f^j_{i,\xi})_\eta \mid i=1,\dots,l;\quad j=0,1,\dots,\deg f_i-1\}
\]
are linearly independent. This already proves that $\cf_\xi$ is a polynomial algebra
freely generated by the $\{f^j_{i,\xi}\}$'s . We have also proved the following 
implication:
\\[.7ex]
{\it if\/ $\bbk\xi+\bbk\eta$ is a $\q^*_{reg}$-plane, then the vectors 
$ \{ (\textsl{d}f^j_{i,\xi})_\eta \mid i=1,\dots,l;\quad j=0,1,\dots,\deg f_i-1\}$
are linearly independent.}

(ii) Now $\codim (\q^*\setminus \Omega)\ge 3$.
Applying Lemma~\ref{auxil} to $V=\q^*$, $\Omega=\q^*_{reg}$, and $v=\xi$, we 
conclude that 
\[
\{\nu\in \q^*_{reg}\mid  
(\textsl{d}f^j_{i,\xi})_\nu\ \text{ are linearly independent}\}
\]
is an open subset
of $\q^*$ whose complement is of codimension $\ge 2$. 
This means, in turn, that  \cite[Theorem~1.1]{ppy} applies to the polynomial subalgebra
$\cf_\xi\subset \cs(\q)$.
Therefore,
we can conclude that  the subalgebra $\cf_\xi$ is algebraically closed in $\cs(\q)$. 

Assume that $\eus K$ is a commutative subalgebra of $\cs(\q)$ containing
$\cf_\xi$. Since $\cf_\xi$ has the maximal possible Krull dimension,
$\cf_\xi\subset \eus K$ is a an algebraic extension. Because $\cf_\xi$ is 
algebraically closed in $\cs(\q)$, we obtain $\cf_\xi=\eus K$.
\end{proof}

\begin{rmk}
The codim--3 property is essential for the maximality of 
$\cf_\xi(\zq)$, see Example~\ref{max_rank}.
\end{rmk}

It would be interesting to find general conditions that guarantee us
that the family of $\xi$-shifts of the free generators of $\zq$ form a
regular sequence in $\cs(\q)$. In the geometric language, this means that
we are interested in the property that the natural morphism
$\q^* \to \spe(\cf_\xi(\zq))\simeq \mathbb A^{b(\q)}$ is flat. It is likely that the assumptions of 
Theorem~\ref{thm:cod3} are sufficient for this. However, we unable to prove this as yet.

\begin{rmk}
One can use deformation arguments for proving flatness. We mention 
an affirmative result for $\sln$, which is obtained by combining work of several
authors. For an arbitrary reductive $\g$,
there is a general procedure of obtaining new commutative subalgebras of $\cs(\g)$ 
as limits of Mishchenko-Fomenko subalgebras $\cf_\xi(\zg)$, where $\xi$ runs
inside a fixed Cartan subalgebra of $\g$, 
see \cite{vitya}.
In particular, for $\g=\sln$, there is a special limit subalgebra that is the associated
graded algebra of the Gelfand-Zetlin subalgebra of $\eus U(\sln)$, see \cite[\S\,6]{vi90}.
In \cite{ovs}, it is proved that the free generators of the latter form a regular sequence in
$\cs(\sln)$. This implies that if $\xi\in(\sln)^*\simeq \sln$ is regular semisimple, 
then the free generators of $\cf_\xi(\eus Z(\sln))$ form a regular sequence.
\end{rmk}

\section{Applications} 
\label{sect:applic}

\subsection{Some Lie algebras with codim--$3$ property}
Here we describe several classes of Lie algebras, where
Theorem~\ref{thm:cod3} applies.

1) \ If $\g$ is reductive, then the assumptions of Theorem~\ref{thm:cod3}
are satisfied. This follows from the classical results of Kostant \cite{ko63}.
Therefore,  for any 
$\xi\in\g_{reg}$, $\cf_\xi(\eus Z(\g))$ is a polynomial algebra, and it is a
maximal commutative subalgebra of $\cs(\g)$.
For the regular semisimple $\xi$, this has already been proved by Tarasov~\cite{tar1}.

2) \ Following \cite{rt}, recall the definition of a (generalised) Takiff Lie algebra
(modelled on $\q$).
The infinite-dimensional $\bbk$-vector space $\q_\infty:=\q\otimes \bbk[\mathsf T]$ has a natural
structure of a Lie algebra such that $[x\otimes \mathsf T^l, y\otimes \mathsf T^k]=
[x,y]\otimes \mathsf T^{l+k}$.
Then $\q_{\ge (n+1)}=\displaystyle\bigoplus_{j\ge n+1} \q\otimes \mathsf T^j$ is an 
ideal of
$\q_\infty$, and $\q_\infty/\q_{\ge (n+1)}$
is a {\it generalised Takiff Lie algebra}, denoted $\q\langle n\rangle$.
If  $\q=\g$ is semisimple, then $\g\langle n\rangle$ satisfies all the assumptions 
of Theorem~\ref{thm:cod3}, see \cite{rt}. For $n=1$, one obtains the semi-direct product
$\g \ltimes\g$. This case was studied by Takiff in 1971.

3) \ Let $e\in\sln$ be a nilpotent element. Set $\q=\z_\sln(e)$.
Then $\ind\q=\rk(\sln)=n{-}1$ \cite{kos} and $\cs(\q)^Q$ is a polynomial algebra 
of Krull dimension $n{-}1$
such that the sum of the degrees of free generators equals $b(\q)$ \cite[Theorem\,4.2]{ppy}.
The second author can prove that here $(\q,\ads)$ have codim--$3$ property.
(This will appear elsewhere.) Thus,
$\z_\sln(e)$ satisfies all the assumptions 
of Theorem~\ref{thm:cod3}.

4) \ Let $\q$ be a $\BZ_2$-contraction of a simple Lie algebra $\g$.
It is known that $\trdeg\zq=\ind\q$ \cite[Lemma\,2.6]{coadj} and
 $(\q,\ads)$ has the codim--$2$ property \cite[Theorem\,3.3]{coadj}.
However, the stronger codim--$3$ property is not always satisfied.
Recall the relevant setup. 

Let $\g=\g_0\oplus\g_1$ be a $\BZ_2$-grading
of $\g$. Then the semi-direct product  $\q=\g_0\ltimes\g_1$ is called  a
$\BZ_2$-{\it contraction\/} of $\g$. Here $\ind\q=\ind\g=\rk\g$, 
hence $b(\q)=b(\g)$. For most $\BZ_2$-gradings, it is proved
that $\zq$ is polynomial
and the sum of degrees of free generators equals $b(\g)$, see \cite[Sect.\,4 \& 5]{coadj}.
It follows that, for such $\BZ_2$-contractions, the commutative subalgebras $\cf_\xi(\zq)$,
$\xi\in\q^*_{reg}$,  are polynomial and 
of maximal dimension. However, these are not always maximal.

\begin{ex}   \label{max_rank}
Let $\g=\g_0\oplus\g_1$ be a $\BZ_2$-grading such that $\g_1$ contains
a Cartan subalgebra of $\g_1$. It is equivalent to that $\dim\g_1=b(\g)$.
Then $\cs(\q)^Q=\cs(\g_1)^{G_0}\simeq \cs(\g)^G$.
(This clearly shows that the sum of degrees of free generators of  $\cs(\q)^Q$
equals $b(\g)$.)
By the assumption, $\g_1$ contains regular elements of $\g$ and, hence, of $\q$.
Let $\xi\in\g_1$ be such an element.
Then $\cf_\xi(\zq)=\cf_\xi(\cs(\g_1)^{G_0})$ is a proper subalgebra of
$\cs(\g_1)$. Indeed, the family of $\xi$-shifts of the generators contains
$b(\g)$ elements, but not all of them are of degree 1. On the other hand,
the subspace $\g_1$ is a commutative Lie subalgebra of $\q$, hence 
$\cs(\g_1)$ is a commutative subalgebra of $\cs(\q)$. (Actually, it is a maximal
commutative subalgebra!) Thus, $\cf_\xi(\zq)$ is a commutative
subalgebra of $\cs(\q)$ of maximal dimension, but not maximal.
\\[.6ex]
Of course, the reason for such a "bad" behaviour is that 
$\codim (\q^*\setminus \q^*_{reg})=2$. This can  also be proved directly using 
invariant-theoretic properties of the $G_0$-module $\g_1$ \cite{kr71}.
\end{ex}

\begin{ex}   \label{ex:good} 
We have verified that
the codim--$3$ property holds for  $\BZ_2$-contractions associated with
the following symmetric pairs $(\g,\g_0)$: \  
$(\mathfrak{sl}_{2n},\spn)$;  $(\slno,\gln)$, $n\ge 2$; $(\son,\mathfrak{so}_{n-1})$; 
$(\GR{E}{6},\GR{F}{4})$;
$(\GR{F}{4},\GR{B}{4})$. However, the complete list is not known yet.
For items~2,3, and 5, it is shown in \cite{coadj} that $\zq$ is polynomial and 
the sum of degrees of the free generators equals $b(\q)$. 
Hence Theorem~\ref{thm:cod3} applies there.
\end{ex}

\begin{rmk} In \cite{jp}, Joseph and Lamprou consider the so-called 
{\it truncated parabolic subalgebras of maximal index  in} $\sln$.
In this case, 
$\zq$ is a polynomial algebra and the equality $\sum \deg f_i=b(\q)$ holds.
Furthermore, they prove the maximality of  $\cf_\xi(\zq)$.
It would be interesting to check whether the codim--$3$ property also holds there.
\end{rmk}

\subsection{Semi-direct products and the codim--$3$ property}
Example~\ref{max_rank} can be put in a more general context.
Suppose $G$ is semisimple  and $V$ is  a finite-dimensional $G$-module. Set 
$m=\max_{\zeta\in V^*}\dim G{\cdot}\zeta$.
Form the semi-direct product $\q=\g\ltimes V$.

\begin{prop}   \label{context}
Suppose that {\sf (a)} \ $S(V)^G=\bbk[V^*]^G$ is a polynomial algebra
and {\sf (b)} \ $m=\dim\g$.
Then $(\q,\ads)$ does not satisfy the codim--$3$ property and
 the commutative subalgebras $\cf_\xi(\zq)$ are not maximal.
\end{prop}\begin{proof}
It follows from assumption~(b) and Ra\"\i s' formula \cite{rais} that $\ind\q=
\dim V-\dim\g$ and therefore $b(\q)=\dim V$. Also, assumption~(b) implies that
$\bbk[\q^*]^{Q}=\bbk[V^*]^G$ \cite[Theorem\,6.4]{p05}. 
Thus, $\zq=\cs(\q)^{Q}=\bbk[\q^*]^{Q}$ is a polynomial 
algebra.  Since $G$ has no rational characters, $\bbk(V^*)^G$ is the quotient field of
$\bbk[V^*]^G$. Hence $\trdeg \bbk[V^*]^G=\ind\q$. Let $d$ be the sum of degrees
of free generators of $\bbk[V^*]^G$. By \cite[Korollar\,6]{knop}, $d\le \dim V$.
Assume that $(\q,\ads)$ has the codim--$3$ property. Then
$d\ge b(\q)=\dim V$ (Theorem~\ref{thm:cod2}). 
Hence $d=b(\q)$ and by Theorem~\ref{thm:cod3},
$\cf_\xi(\zq)$ is a maximal commutative subalgebra of $\cs(\q)$ for any $\xi\in\q^*_{reg}$.
Since $\zq$ is a subalgebra of $\cs(V)$, $\cf_\xi(\zq)$ is a  subalgebra of $\cs(V)$, too.
Furthermore, $\cf_\xi(\zq)$ is generated by $\dim V$ elements, and not all of them are of
degree 1. Thus,  $\cf_\xi(\zq)$ is a proper subalgebra of $\cs(V)$,
and the latter is a (maximal) commutative subalgebra of $\cs(\q)$.
This contradiction shows that the codim--$3$ property cannot be satisfied for 
$(\q,\ads)$. The above argument also proves the second assertion.
\end{proof}%
\begin{rmk} Set $V^*_{sing}=\{\nu\in V^* \mid \dim G{\cdot}\nu < m\}$. 
(This closed subset plays an important r\^ole in theory developed in \cite{knop}.) 
It is easily seen that if $m=\dim G$ and
$\codim V^*_{sing}\ge n$, then $\codim \q^*\setminus \q^*_{reg}\ge n$.
Hence, under the assumptions of Proposition~\ref{context}, we have $\codim V^*_{sing}\le 2$,
and according to \cite[Korollar\,2]{knop},  $\codim V^*_{sing}= 2$ if and only if
$d=b(\q)$.
\end{rmk}

\appendix   
\section{Some results on skew-symmetric bilinear forms}   
\label{app}
\setcounter{equation}{0}

\noindent
Here we present some general facts concerning  skew-symmetric bilinear forms
that are needed for the proof of Theorem~\ref{compl}. 
All these results are extracted from \cite{bol2}, but we present them in a more systematic
and algebraic form.

Let $\eus P$
be a two-dimensional linear space of (possibly degenerate) skew-symmetric bilinear
forms on a finite-dimensional vector space $V$. Set $m=\max_{A\in \eus P }\rk A$, and let 
$\eus P_{reg}\subset \eus P$ be the set of all forms of rank $m$.
For each $A\in \eus P$, let $\ker A\subset V$ be the kernel of $A$. 
Our main object of interest is the subspace 
$L:=\sum_{A\in \eus P_{reg}} \ker A$.  

\begin{lm}              \label{open}
For any nonempty open subset\/ $\Omega\subset \eus P_{reg}$, 
we have $\sum_{A\in \Omega} \ker A=L$.
\end{lm}
\begin{proof} 
Set $r=\dim V-m$ and $M=\sum_{A\in \Omega} \ker A \subset L$.
Take any $C\in \eus P_{reg}\setminus \Omega$.  Then $\ker C$ is a point of  
the Grassmannian ${\rm Gr}_r(V)$. Because  $\eus P$ is irreducible, 
$\overline{\Omega}=\eus P$ and 
there is a curve $\varkappa :\bbk^\times \to \Omega$
such that $\lim_{t\to 0} \varkappa(t) = C$. Hence 
\[
\lim_{t\to 0}(\ker\varkappa(t))=\ker C,
\]
where the last limit is taken in ${\rm Gr}_r(V)$. Since
 $\ker\varkappa(t)\in {\rm Gr}_r(M)$ for $t\ne 0$ and 
 ${\rm Gr}_r(M)$ is closed in ${\rm Gr}_r(V)$, we obtain 
$\ker C\subset M$. Thus, $M=L$.
\end{proof}

\noindent
For $A\in \eus P$, let $\hat A$ denote the corresponding linear
map from $V$ to $V^*$. Then $\ker A=\ker\hat A$.

\begin{lm}           \label{image}
For all $A,B\in \eus P\setminus\{ 0\}$, we have $\hat A(L)=\hat B(L)$. 
\end{lm}
\begin{proof} Clearly, we may assume that $A$ and $B$ are linearly independent.
By virtue of Lemma~\ref{open}, 
$L$ is spanned by some $L_{a,b}:=\ker(a A+b B)$ with 
$ab\ne 0$. 
Since $(a\hat A+b\hat B)(L_{a,b})=0$, we obtain 
$(a\hat A)(L_{a,b})=(b\hat B)(L_{a,b})$ and hence $\hat A(L_{a,b})=\hat B(L_{a,b})$. 
The result follows. 
\end{proof}

\noindent
For $A\in \eus P\setminus\{0\}$, let $\widetilde{L}\subset V$ denote
the annihilator of $\hat A(L)\subset V^*$.  By Lemma~\ref{image},
$\widetilde{L}$ does not depend on the choice of $A$.
Note also that 
$\widetilde{L}=\{v\in V\mid A(v,L)=0\}$.
Since $\ker A\subset \widetilde{L}$ for each nonzero $A$, 
$L$ is a subspace of $\widetilde{L}$. 

\begin{lm}  \label{phi}
Suppose that $B\in \eus P$ and $A\in \eus P_{reg}$. Then 
\begin{itemize}
\item[\sf (i)] \ $\hat B(\widetilde{L})\subset \hat A(\widetilde{L})$;
\item[\sf (ii)] \ Associated with $A$ and $B$, there is a natural linear operator
$\Phi_{A,B}=\Phi: \widetilde{L}/ L \to \widetilde{L}/ L$.
\end{itemize}
\end{lm}
\begin{proof}
(i) Let $M_A$ and $M_B$ be the the annihilators of 
$\hat A(\widetilde{L})$ and $\hat B(\widetilde{L})$, respectively.
Since  
$M_A=\ker A+L=L$ and $M_B=\ker B+L$, we obtain $M_A\subset M_B$.

(ii) Take any $v\in\widetilde{L}$. 
Since $\hat B(\widetilde{L})\subset \hat A(\widetilde{L})$, 
where is $w\in\widetilde{L}$ such that $\hat A(w)=\hat B(v)$. 
Letting $\Phi(v+L):=w+L$, we have to check that there is no ambiguity
in this. To this end, assume  that 
$\hat A(w')\in \hat B(v+L)=\hat A(w)+\hat B(L)$.
Since $\hat B(L)=\hat A(L)$, we obtain $\hat A(w'-w)\in \hat A(L)$.
Hence $w-w'\in L+\ker A=L$. Thus, given 
$\bar v=v+L\in\widetilde{L}/L$, there is a
unique $\bar w=w+L\in\widetilde{L}/L$ such that $\hat B(\bar v)=\hat A(\bar w)$. 
The claim follows.
\end{proof}
\begin{prop}
\label{com1}
If $\eus P_{reg}=\eus P\setminus\{0\}$, then $L=\widetilde{L}$; in other words,
$L$ is a maximal isotropic subspace of\/ $V$ with respect  to any nonzero
$A\in \eus P$.
\end{prop} 
\begin{proof} Take linearly independent 
$A$ and $B$, as in Lemma~\ref{phi}.
We use  the operator  $\Phi: \widetilde{L}/L\to \widetilde{L}/L$ 
introduced in Lemma~\ref{phi}(ii). Since $\bbk$ is algebraically closed, 
$\widetilde{L}/L=\{0\}$ if and only if all eigenvectors of $\Phi$ are zero. 
Assume that $v+L\in\widetilde{L}/L$ is a  $\lb$-eigenvector of $\Phi$. 
Then expanding the definition of $\Phi$ yields  
$(\hat B-\lb \hat A)v \in \hat A(L)$. Since $\hat A(L)=(\hat B-\lb \hat A)(L)$
by Lemma~\ref{image}, we get  
$(\hat B-\lb \hat A)(v)\in (\hat B-\lb \hat A)(L)$ and, hence,  $v\in L+\ker(B-\lb A)$.
If $v\not\in L$, then 
$\ker (B-\lb A)\not\subset L$ and therefore
$(B-\lambda A)\not\in \eus P_{reg}$. A contradiction!
\end{proof}

\end{document}